\documentclass[12pt]{article}
\usepackage{amssymb}
\begin{document}
\title{Properly discontinuous actions on bounded domains}
\author{Bo-Yong Chen\footnote{Address: Department of Mathematics, Tongji University, Shanghai
200092, P. R. China. E-mail: boychen@mail.tongji.edu.cn}}
\date{February 14, 2008}
\maketitle

\begin{abstract}
We give sufficient conditions for the quotient of a free, properly
discontinuous action on a bounded domain of holomorphy to be a Stein
manifold in terms of Poincar\'e series or limit sets for orbits. An
immediate consequence is that the quotient of any cyclic, free,
properly discontinuous action on the unit ball or the bidisc is
Stein.
\end{abstract}
\section{Introduction}

\ \ \ \ The uniformization theorem tells us that one can study
algebraic properties of compact Riemann surfaces with genus greater
than one through the function theory of the unit disc (e.g., the
Poincar\'e series). Such an idea may be extended without any
difficulty to
 high dimensional compact complex manifolds covered by bounded
domains in ${\bf C}^n$ (cf. [13]). However, those non-compact
complex manifolds uniformized by bounded domains, which are much
more than compact ones, are largely ignored, expect the special case
when they can be compactified. Function theory of such manifolds
should be of great interests. In this paper, we shall first prove
some results along this line.

\bigskip

\textbf{Theorem 1.} \emph{Let $\Gamma$ be a free, properly
discontinuous group of automorphisms of the unit ball ${\bf
B}\subset {\bf C}^n$. Suppose there exists a point $a\in {\bf B}$
such that}
\begin{equation}
\sum_{\gamma\in \Gamma} (1-|\gamma(a)|)<\infty.
\end{equation}
\emph{Then ${\bf B}/\Gamma$ is a Stein manifold which possesses a
negative and $C^\infty$ strictly plurisubharmonic function.}

\bigskip

A standard application of the $L^2-$estimate for the
$\bar{\partial}-$operator on complete K\"ahler manifolds (cf. [1],
[3]) yields

\bigskip

\textbf{Corollary 1.} \emph{The Bergman metric exists on ${\bf
B}/\Gamma$}.

\bigskip

\textbf{Remark.} One might expect that condition (1) is not
necessary for the existence of the Bergman metric on ${\bf
B}/\Gamma$. However, Mumford [9] has constructed a two-dimensional
compact ball quotient without non-zero holomorphic $(2,0)-$forms.

\bigskip

In the classification of open Riemann surfaces, one calls an open
Riemann surface hyperbolic if the Green function exists. We have the
following new characterization of hyperbolic Riemann surfaces.

\bigskip

\textbf{Corollary 2.} \emph{An open Riemann surface is hyperbolic
iff it carries a negative, $C^\infty$ strictly subharmonic
function}.

\bigskip

\emph{Proof.} We only need to verify the "only if" part. By the
uniformization theorem we can write a hyperbolic Riemann surface as
$\Delta/\Gamma$ where $\Gamma$ acts freely and properly
discontinuously on the unit disc $\Delta$. According to a classical
theorem of Myrberg (cf. Theorem XI. 13. in [15]),  $\Gamma$
satisfies condition (1) and the assertion follows immediately from
Theorem 1.

\bigskip

In $\S\,4$ we will show that every cyclic, properly discontinuous
action of the unit disc $\Delta$ satisfies condition (1), while for
the case of the unit ball ${\bf B}\subset {\bf C}^n$, $n>1$ most
cyclic, properly discontinuous actions of ${\bf B}$ satisfies
condition (1) except when generator is certain parabolic elements.

\bigskip

The method for proving Theorem 1 also works for those bounded
domains admitting a negative \emph{Strictly} plurisubharmonic
exhaustion function whose modulus is asymptotic to the boundary
distance (e.g. strongly pseudoconvex domains), however it is not
sufficient for more general pseudoconvex domains (e.g., bounded
symmetric domains). Thus a variation of this approach is needed.

\bigskip

In order to present a general result, we need some preparation. Let
$M$ be a complex manifold of dimension $n$. The \emph{pluricomplex
Green function} with pole $w$ is defined by
$$
g_M(z,w)=\sup\{u(z)\}
$$
where the supremum is taken over all negative plurisubharmonic
functions $u$ on $M$ satisfying $u(z)\le \log|z-w|+O(1)$ in a
suitable coordinate neighborhood around $w$. If $F:M\rightarrow M'$
is a holomorphic map, then $g_M(z,w)\ge g_{M'}(F(z),F(w))$ and
equality holds when $F$ is biholomorphic. We say that $g_M(\cdot, w
)$ exists if
$$
|g_M(z, w )-\log|z-w||=O(1)
$$
holds in a coordinate neighborhood of $w$ and $-g_M(z,w)$ is bounded
outside this neighborhood (i.e., $g_M(z,w)$ has no other pole than
$w$). $M$ is called \emph{pluricomplex hyperbolic} if $g_M(\cdot, w
)$ exists for every $w\in M$. Typical examples involve bounded
domains in ${\bf C}^n$ and their biholomorphic images. It is also
known that the Bergman metric exists on pluricomplex hyperbolic
Stein manifolds (cf. [2]).

\bigskip

\textbf{Theorem 2.} \emph{Let $M$ be a pluricomplex hyperbolic Stein
manifold. Suppose $\Gamma$ is a free, properly discontinuous group
of automorphisms of $M$ so that the inequality }
\begin{equation}
\sum_{\gamma\in \Gamma,\gamma\neq I} g_M(\gamma(z),z)>-\infty
\end{equation}
\emph{for all $z\in M$. Then $M/\Gamma$ is also a pluricomplex
hyperbolic Stein manifold. }

\bigskip

\textbf{Theorem 3.} \emph{The condition of Theorem 2 is satisfied in
the following cases:}

1) \emph{$M:$ bounded convex domain in ${\bf C}^n$, $\Gamma:$ free,
 properly discontinuous group so that the inequality}
$$
\sum_{\gamma\in \Gamma}\delta_M(\gamma(a))<\infty
$$
\emph{holds for some $a\in M$. Here $\delta_M$ denotes the Euclidean
boundary distance. }

2) \emph{$M:$ the Teichm\"uller space of Riemann surfaces with genus
$g$ and $n-$punctures, $\Gamma:$ free,
 properly discontinuous group so that the inequality}
$$
\sum_{\gamma\in \Gamma} e^{-2d_T(\gamma(a),a)}<\infty
$$
\emph{holds for some $a\in M$. Here $d_T$ is the Teichm\"uller
distance.}

 \bigskip

 Next, we study how the limit set of a
 properly  discontinuous action $\Gamma$ on a bounded domain $\Omega$ in ${\bf C}^n$ infects the geometry of $\Omega$ or $\Omega/\Gamma$.
A point $w\in
\partial \Omega$ is called \emph{a limit point} of $\Gamma$ if there
exists a compact subset $K$, $\{z_k\}\subset \Omega$ with
$z_k\rightarrow w$ and $\{\gamma_k\}\subset \Gamma$ such that
$\gamma_k(z_k)\in K$.  Let $L_\Gamma$ be the set of all limit points
of $\Gamma$. For a fundamental domain $D$ of $\Gamma$ and $a\in D$,
we define by $\Gamma_a$ the orbit at $a$ and call ${\rm C}\Gamma_a$
the set of \emph{cluster points} of $\Gamma_a$. It is obvious that
$\cup_{a\in D}{\rm C}\Gamma_a\subset L_\Gamma$.
  Let $K_\Omega$ denote the Bergman kernel of $\Omega$.
 Define
 $$
 \Lambda_\Omega=\left\{w\in \partial \Omega:{\lim\sup}_{z\rightarrow
 w}K_\Omega(z)=\infty\right\}
 $$
 where $\partial \Omega$ is the boundary of $\Omega$. We call
 $\Omega$ \emph{an $L^2_h-$domain of holomorphy} if there
 are no domains $\Omega_1,\Omega_2$ in ${\bf C}^n$ with $\Omega_2\subset \Omega\cap \Omega_1,\,({\bf C}^n\backslash\Omega)\cap \Omega_1\neq \emptyset$
so that for any $L^2-$holomorphic function $f$ in $\Omega$ there
exists a holomorphic function $F$ in $\Omega_1$ with $F=f$ in
$\Omega_2$. Clearly, every $L^2_h-$domain of holomorphy must be a
domain of holomorphy. It is also known that $\Omega$ is an
$L^2_h-$domain of holomorphy if and only if $\Lambda_\Omega=\partial
\Omega$ (cf. [10]).

 \bigskip

 \textbf{Theorem 4.} \emph{$L_\Gamma\subset \Lambda_\Omega$. In particular,
 if $L_\Gamma$ is dense in $\partial \Omega$ $($in the Euclidean topology$)$, then $\Omega$ is an $L^2_h-$domain of holomorphy.
 }

 \bigskip

Clearly, if $\Gamma$ has a compact fundamental domain, then
$L_\Gamma=\partial \Omega$. Thus we get a generalization of a
well-known theorem of Siegel on Steinness of a bounded domain  when
$\Omega$ has a compact quotient (cf. [13]).
 It is expected that $L_\Gamma$ lies
dense in $\partial \Omega$ when the fundamental domain is
non-compact but the part intersecting $\partial \Omega$ is
sufficiently "small" (eg., finite sets).

\bigskip

The case when the limit set is relatively "small" is of independent
interest, we have the following

\bigskip

\textbf{Theorem 5.} \emph{Let $\Omega$  be a bounded Stein domain in
${\bf C}^n$ and $\Gamma$ a free, properly discontinuous action on
$\Omega$. Suppose there exists a negative continuous
plurisubharmonic function $\psi$ on $\Omega$ so that $\psi>-\infty$
on $\Omega$ and $\psi(z)\rightarrow -\infty$ as $z\rightarrow
\cup_{a\in D} {\rm C}\Gamma_a$ where $D$ is a fundamental domain of
$\Gamma$. Then $\Omega/\Gamma$ is a Stein manifold which carries a
negative, $C^\infty$ strictly plurisubharmonic function. In
particular, the Bergman metric exists on $\Omega/\Gamma$.}

\bigskip

Clearly, the condition of Theorem 5 is satisfied when $\cup_{a\in D}
{\rm C}\Gamma_a$ is complete pluripolar, i.e., it is precisely the
set for certain plurisubharmonic function in ${\bf C}^n$ taking
value $-\infty$ (e.g., finite sets). Based on the following results,
it seems to be true that free, properly discontinuous and cyclic
actions of any bounded Stein domain $\Omega$ would satisfy this
condition.

\bigskip

\textbf{Theorem 6.} \emph{For every infinite, cyclic, free and
properly discontinuous action $\Gamma$ of ${\bf B}$, the set
$\cup_{a\in D} {\rm C}\Gamma_a$ has at most two elements. In
particular, the quotient is a Stein manifold}.

\bigskip

\textbf{Theorem 7.} \emph{Every infinite, cyclic, free and properly
discontinuous action  of the bidisc $\Delta^2$ satisfies the
condition of theorem 5. In particular, the quotient is a Stein
manifold}.

\bigskip

\textbf{Acknowledgement.} After the manuscript was completed, the
author became aware of the existence of [6] where the Steinness  of
quotients of the unit ball for cyclic, free and properly
discontinuous actions is proved (actually they are even
biholomorphic to bounded Stein domains!). He thanks Professor Takeo
Ohsawa for pointing out this reference and valuable suggestions.

\section{Preliminaries}

\ \ \ \ Let $M$ be a complex manifold and ${\rm Aut}(M)$ the group
of automorphisms of $M$. A subgroup $\Gamma$ of ${\rm Aut}(M)$ is
called \emph{properly discontinuous} if for any compact subsets
$K_1,\,K_2$ of $M$ there are only a finite number of elements
$\gamma\in \Gamma$ such that $\gamma(K_1)\cap K_2\neq \emptyset$.
Two points $z,z'$ in $M$ are called \emph{equivalent} if
$z'=\gamma(z)$ for some $\gamma\in \Gamma$. A set $D\subset M$ is
called \emph{a fundamental domain} of $M$ if any two points of $D$
are not equivalent and any point $z\in M$ has its equivalent in $D$.
In case when $M$ is a bounded domain in ${\bf C}^n$, $\Gamma$ is
properly discontinuous iff  there exists one point $a\in M$ such
that equivalents of $a$ have no cluster points in $M$ (cf. Chapter
13 in [13]).

\bigskip

\textbf{Proposition 1.} \emph{Let $\Omega$ be a bounded domain  in
${\bf C}^n$ and $\Gamma$ a countable subgroup of ${\rm
Aut}(\Omega)$. Let $\tau$ be a positive continuous function on
$\Omega$. Suppose that there exists $a\in \Omega$ such that}
\begin{equation}
\sum_{\gamma\in \Gamma}\tau(\gamma(a))<\infty,
\end{equation}
\emph{then $\Gamma$ acts properly discontinuously on $\Omega$. }

\bigskip

\emph{Proof.} Suppose $\Gamma$ is not properly discontinuous, then
the orbit $\Gamma_a$ would have at least one cluster point $b$. By
the continuity of $\tau$ there are infinite $\gamma\in \Gamma$ so
that
$$
\tau(\gamma(a))\ge \frac12 \tau(b)>0,
$$
contradicts with $(3)$.

\bigskip

For a bounded domain $\Omega$ and a properly discontinuous action
$\Gamma$ on $\Omega$, I do  not known whether the limit set
$L_\Gamma$ coincides with $\cup_{a\in D}{\rm C}\Gamma_a$ ($D:$  a
fundamental domain), except for the following special case.  For a
point $p\in
\partial \Omega$, if there is a function $\psi$, continuous on the
closure of $\Omega$, plurisubharmonic in $\Omega$, such that
$\psi(p)=0$ and $\psi(z)<0$ for every $z\in
\overline{\Omega}\backslash\{p\}$, then we call $\psi$ \emph{a
barrier at $p$}. A domain $\Omega$ is called \emph{B-regular} if
there exists a barrier at every boundary point.

\bigskip

\textbf{Proposition 2.} \emph{Let $\Omega$ be a bounded B-regular
domain in ${\bf C}^n$ and $\Gamma$ a properly discontinuous action
on $\Omega$. Then ${\rm C}\Gamma_a=L_\Gamma$ for all $a\in D$.}

\bigskip

\emph{Proof.} Clearly, ${\rm C}\Gamma_a\subset L_\Gamma$ for every
$a\in D$. On the other hand, for every $p\in L_\Gamma$ there is  a
compact subset $K$ of $\Omega$, $\{\gamma_k\}\subset \Gamma$ and
$\Omega \ni p_k\rightarrow p$ satisfying $\gamma_k(p_k)\in K$. Since
$\Omega$ is bounded, $\{\gamma_k^{-1}\}$ is equicontinuous on every
compact subset of $\Omega$ so that it has at least one limit point,
say $\beta$. After taking a subsequence, we may assume
$\beta(\gamma_k(p_k))\rightarrow p$ as $k\rightarrow \infty$.

Let $\psi$ be a barrier at $p$ and put $\varphi=\psi\circ \beta$.
Then $\varphi$ is plurisubharmonic in $\Omega$, $\varphi\le 0$ on
$\Omega$ and
$$
\max_K \varphi=0.
$$
By the maximum principle, $\varphi=0$ on $\Omega$, this means
$\gamma^{-1}_k(a)\rightarrow p$ for every $a\in D$ so that $p\in
{\rm C}\Gamma_a$.

\bigskip

\textbf{Proposition 3.} \emph{Let $\Omega$ be a bounded domain which
has a barrier at some $p\in \partial \Omega$. Let $q\in \Omega$ and
$\{\gamma_k\}\subset {\rm Aut}(\Omega)$ satisfying
$\gamma_k(q)\rightarrow p$ as $k\rightarrow \infty$. Then
$\gamma_k(z)\rightarrow p$ uniformly on compact subsets of
$\Omega$.}

\bigskip

\emph{Proof.} The argument is quite similar as above. Let $\gamma$
be any limit point of $\{\gamma_k\}$ in the topology of uniform
convergence on compact subsets of $\Omega$. Then $\gamma(q)=p$. Put
$\varphi=\psi\circ \gamma$ ($\psi:$ the barrier at $p$). Similarly,
we have $\varphi=0$ on $\Omega$, i.e., $\gamma(z)=p$ for all $z\in
\Omega$. This means that the point $p$ is the only limit of
$\{\gamma_k(z)\}$ for every $z$.

\section{Proof of Theorem 1}

\ \ \ \ 3.1. Fix $a\in {\bf B}$. Let $P_0=0$ and
$$
P_a (z)=\frac{\langle z,a\rangle}{|a|^2}a \ \ \ {\rm if\ } a\neq 0,\
\ \ Q_a=I-P_a
$$
where $\langle z,w\rangle=\sum_{j=1}^n z_j \bar{w}_j$. Put
$s_a=(1-|a|^2)^{1/2}$ and define
$$
\varphi_a(z)=\frac{a-P_a(z)-s_aQ_a(z)}{1-\langle z,a\rangle}.
$$
We have the following well-known properties: $\varphi_a(0)=a$,
$\varphi_a(a)=0$ and the identity
\begin{equation}
1-|\varphi_a(z)|^2=\frac{(1-|a|^2)(1-|z|^2)}{|1-\langle
z,a\rangle|^2}
\end{equation}
holds for every $z\in \overline{{\bf B}}$ (cf. Chapter 2 of [12]).
Let ${\rm U}(n)$ denote the group of the unitary transformations in
${\bf C}^n$. Then
$$
{\rm Aut}({\bf B})=\{U\varphi_a:a\in {\bf B},\ U\in {\rm U}(n)\}.
$$

\bigskip

3.2. Let $\Gamma$ be a properly discontinuous group of ${\rm
Aut}({\bf B})$ and let $a_\gamma=\gamma(a)$, $z_\gamma=\gamma(0)$.
 Recall that the Bergman distance
from $0$ to a point $z\in {\bf B}$ is given by
$$
d_{\bf B}(0,z)=\frac{\sqrt{n+1}}2 \log \frac{1+|z|}{1-|z|}.
$$
Since $d_{\bf B}(a,0)=d_{\bf B}(a_\gamma,z_\gamma)=d_{\bf
B}(a_\gamma,\varphi_{z_\gamma}^{-1}(0))=d_{\bf
B}(\varphi_{z_\gamma}(a_\gamma),0)$, we infer from (4) that
\begin{eqnarray*}
&& 1-|a|^2=1-|\varphi_{z_\gamma}(a_\gamma)|^2  =
\frac{(1-|z_\gamma|^2)(1-|a_\gamma|^2)}{|1-\langle a_\gamma,z_\gamma\rangle|^2}\\
& \le &
4\min\left\{\frac{(1-|a_\gamma|)(1-|z_\gamma|)}{(1-|z_\gamma|)^2},\frac{(1-|a_\gamma|)(1-|z_\gamma|)}{(1-|a_\gamma|)^2}\right\},
\end{eqnarray*}
so that
$$
\frac{(1-|a|^2)(1-|z_\gamma|)}4 \le 1-|a_\gamma|\le
\frac{4(1-|z_\gamma|)}{1-|a|^2}.
$$
From this we conclude that if condition (1) holds at one point in
${\bf B}$ then it holds for every point in ${\bf B}$.

\bigskip

3.3.  Put
$$
u(z)=\sum_{\gamma\in \Gamma} (|\gamma^{-1}(z)|^2-1),\ \ \ \ \ z\in
{\bf B}.
$$
Clearly, we have $u(\gamma(z))=u(z)$ for every $\gamma\in \Gamma$.
 Fix $0<r<1$ and a point $z\in {\bf B}$ with $|z|\le r$. Put
$z_\gamma=\gamma(0)$. By 3.2, we may assume
$$
\sum_{\gamma\in \Gamma} (1-|z_\gamma|)<\infty.
$$
Similar as above,  since $d_{\bf B}(\gamma^{-1}(z),0)=d_{\bf
B}(z,z_\gamma)=d_{\bf B}(\varphi_{z_\gamma}(z),0)$, we have
$$
 1-|\gamma^{-1}(z)|^2=1-|\varphi_{z_\gamma}(z)|^2  =
\frac{(1-|z_\gamma|^2)(1-|z|^2)}{|1-\langle z,z_\gamma\rangle|^2}\le
\frac{4(1-|z_\gamma|)}{1-|z|},
$$
which implies the convergence of $u$. Consequently, $u$ descends to
a negative $C^\infty$ strictly plurisubharmonic function on ${\bf
B}/\Gamma$.

\bigskip

3.4. \emph{Steiness.} We need the following two important results:

\bigskip

\textbf{Ellencwajg Theorem.} (cf. [5]) \emph{If $X$ is a complex
manifold admitting a continuous strictly plurisubharmonic function,
then every locally Stein relatively compact open sets $\Omega$ is
Stein.}

\bigskip

\textbf{Docquier-Grauert Theorem.} (cf. [4]) \emph{The union of an
increasing 1-parameter family of Stein manifolds is Stein.}

\bigskip

Put ${\bf B}_t=\{z\in {\bf C}^n:|z|<t\}$, $0<t<1$. Then ${\bf
B}=\bigcup_{t\in (0,1)}{\bf B}_t$. Let $\pi:{\bf B}\rightarrow {\bf
B}/\Gamma$ be the canonical projection and define $V_t=\pi({\bf
B}_t)$. Clearly, we have ${\bf B}/\Gamma=\bigcup_{t\in (0,1)}V_t$
with $V_t\subset V_{t'}$ for arbitrary $t<t'$. Since $\pi$ is an
open mapping, $V_t$ is open and relatively compact in ${\bf
B}/\Gamma$. Note that $V_t$ is locally Stein since $\pi$ is locally
biholomorphic, hence is Stein according to the Ellencwajg Theorem
since there is a $C^\infty$ strictly plurisubharmonic function on
${\bf B}/\Gamma$. It follows immediately from the Docquier-Grauert
Theorem that ${\bf B}/\Gamma$ is Stein.

\bigskip

3.5. Finally, let us remark the following classical fact: every
free, properly discontinuous action $\Gamma$ of ${\bf B}$ satisfies
$$
\sum_{\gamma\in \Gamma}(1-|\gamma(z)|)^s<\infty
$$
for all $s>n$ and $z\in {\bf B}$; in case when ${\bf B}/\Gamma$ is
compact or of finite volume w.r.t. the Bergman metric, one has
$$
\sum_{\gamma\in \Gamma}(1-|\gamma(z)|)^n=\infty, \ \ \ z\in {\bf B}
$$
(compare Theorem XI. 8, 10 in [15]). It would be interesting to know
how the convergence of the above Poincar\'e series infects the
function theory of ${\bf B}/\Gamma$ for $1<s<n$.
\section{Cyclic actions on ${\bf B}$}

\ \ \ \ In this section we will check condition (1) for a large
class of cyclic, properly discontinuous actions on the unit ball in
${\bf C}^n$.

 $(i)$ $n>1$. Consider the Cayley transform
$$
\Phi:{\bf B}\rightarrow {\bf H}:=\{z\in {\bf C}^n:{\rm
Im\,}z^1>|z'|^2\}
$$
$$
z\rightarrow \left(i\frac{1+z^1}{1-z^1},\frac{iz'}{1-z^1}\right)
$$
where $z'=(z^2,\cdots,z^n)$ and $|z'|^2=|z^2|^2+\cdots+|z^n|^2$.
Clearly, we have
$$
\Phi^{-1}(z)=\left(\frac{z^1-i}{z^1+i},\frac{2z'}{z^1+i}\right).
$$
Note that ${\bf H}$ has a transitive group of automorphisms
generated by the following two type of
transforms:\\
1) $z\rightarrow (z^1+c+2i\langle z',a'\rangle+i|a'|^2,z'+a')$,
where $c\in {\bf
R} $ and $a'\in {\bf C}^{n-1}$;\\
2) $z\rightarrow (r^2z^1,rz'U')$ where $0<r<\infty$ and $U'\in {\rm
U}(n-1)$.\\
A general element can be written as
$$
\tilde{\gamma}(z)=(r^2z^1+c+2i\langle
rz'U',a'\rangle+i|a'|^2,rz'U'+a').
$$
Consider the cyclic group $\Gamma=\{\gamma^k:k\in {\bf Z}\}$ where
$\gamma=\Phi^{-1}\circ\tilde{\gamma}\circ\Phi$.

$(i)_1$: $r\neq 1$. In this case, we claim firstly that $a'=0$. Here
we only deal with the case $r<1$ (the case $r>1$ is similar).
 Simple computations show that
$$
\tilde{\gamma}^k(i,0')\rightarrow ((1-r^2)^{-1}[c+i|a'|^2+2i\langle
a'(I-rU')^{-1},a'\rangle],a'(I-rU')^{-1}),\ \ \ (k\rightarrow
+\infty).
$$
If $a'\neq 0$, we would have
\begin{eqnarray*}
&& {\rm Im}\left\{\frac1{1-r^2}[c+i|a'|^2+2i\langle a'(I-rU')^{-1},a'\rangle]\right\}-|a'(I-rU')^{-1}|^2\\
& =& \frac1{1-r^2}[|b'(I-rU')|^2+2{\rm
Re}\langle b',b'(I-rU')\rangle-(1-r^2)|b'|^2]\ \ \ (b':=a'(I-rU')^{-1})\\
& = & \frac1{1-r^2}[(2+2r^2)|b'|^2-4{\rm Re}\langle b',b'rU'\rangle]\\
& \ge & \frac{2(1-r)}{1+r}|b'|^2>0,
\end{eqnarray*}
thanks to Schwarz's\ inequality. But this means that
$\{\gamma_k(0)\}$ has a cluster point in ${\bf B}$, absurd. Hence we
are forced to have
$$
\tilde{\gamma}(z)=(r^2z^1+c,rz'U')
$$
and $\tilde{\gamma}^k(i,0')=(r^{2k}i+c(1-r^{2k})(1-r^2)^{-1},0')$.
Hence
\begin{eqnarray*}
&& \sum_{k\in {\bf Z}}(1-|\gamma^k(0)|) \le \sum_{k\in {\bf
Z}}(1-|\gamma^k(0)|^2)\\
 & = & \sum_{k\in {\bf
Z}}\left(1-\left|\frac{(r^{2k}-1)i+c(1-r^{2k})(1-r^2)^{-1}}{(r^{2k}+1)i+c(1-r^{2k})(1-r^2)^{-1}}\right|^2\right)\\
& = & \sum_{k\in {\bf Z}}
\frac{4r^{2k}}{|(r^{2k}+1)i+c(1-r^{2k})(1-r^2)^{-1}|^2}<\infty.
\end{eqnarray*}

$(i)_2:$ $r=1$. In this case, we have
$$
\tilde{\gamma}(z)=(z^1+c+2i\langle z'U',a'\rangle+i|a'|^2,z'U'+a').
$$
A direct computation yields
$$
\tilde{\gamma}^k(i,0')=(i+k(c+i|a'|^2)+2i\langle
a'T_k,a'\rangle,a'S_{k})
$$
where $S_k=I'+U'+\cdots+U'^{k-1}$ and $T_k=(S_1+\cdots+S_{k-1})U'$
($I':$ the identity matrix in ${\bf C}^{n-1}$). Hence
$$
\gamma^k(0)=\left(\frac{kc+ik|a'|^2+2i\langle
a'T_k,a'\rangle}{kc+i[2+k|a'|^2]+2i\langle
a'T_k,a'\rangle},\frac{2a'S_{k}}{kc+i[2+k|a'|^2]+2i\langle
a'T_k,a'\rangle}\right)
$$

$(i)_{21}:$ $a'=0'$. As $\Gamma$ acts properly discontinuously, we
conclude that $c\neq 0$. Therefore,
$$
\sum_{k\in {\bf Z}}(1-|\gamma^k(0)|^2)=\sum_{k\in {\bf
Z}}\left(1-\frac{k^2c^2}{|kc+2i|^2}\right)=\sum_{k\in {\bf
Z}}\frac4{k^2c^2+4}<\infty.
$$

$(i)_{22}:$ $a'\neq 0'$ and $U'=I'$. Then we have
$$
\gamma^k(0)=\left(\frac{kc+ik^2|a'|^2}{kc+i[2+k^2|a'|^2]},\frac{2ka'}{kc+i[2+k^2|a'|^2]}\right).
$$
It is easy to verify
$$
\sum_{k\in {\bf Z}}(1-|\gamma^k(0)|^2)<\infty.
$$

$(i)_{23}:$  $I'-U'$ is invertible and $\langle a'(I'-U')^{-1}U',a'
\rangle\neq -|a'|^2/2$. Simple computations show
$$
S_k=(I'-U'^k)(I'-U')^{-1},\ \ \ T_k=(kI'-S_k)(I'-U')^{-1}U'.
$$
Hence we may write
$$
\gamma^k(0)=\left(\frac{kc+\lambda k+\mu_k}{kc+2i+\lambda
k+\mu_k},\frac{2a'S_{k}}{kc+2i+\lambda k+\mu_k}\right)
$$
where $\lambda\in {\bf C}$ with ${\rm Im\,}\lambda> 0$ and $\mu_k$
are uniformly bounded. Thus
$$
\sum_{k\in {\bf Z}}(1-|\gamma^k(0)|^2)=\sum_{k\in {\bf
Z}}\frac{4k{\rm Im\,}\lambda+O(1)}{k^2|c+\lambda|^2+O(k)}=\infty.
$$

$(i)_{24}:$ $a'\neq 0'$, $I'-U'$ is invertible and $\langle
a'(I'-U')^{-1}U',a' \rangle=-|a'|^2/2$. In this case, the
convergence property holds by a similar argument as
 $(i)_{21}$.

 \bigskip

Generally, one can choose suitable unitary matrix $T\in {\rm
U}(n-1)$ so that $TU'T^{-1}$ becomes a diagonal unitary matrix of
form $[I_{m-1},U_{n-m}]$ where all diagonal elements of $U_{n-m}$
are not 1. Replacing $a'$ by $a'T$, it suffices to consider the case
$U=[I_{m-1},U_{n-m}]$.
 Thus we can reduce the argument to the above cases.

 \bigskip

\emph{Example.} In case $n=2$, $\Gamma$ does not satisfy condition
(1) when
$$
\tilde{\gamma}(z^1,z^2)=(z^1+c+2ie^{i\theta}z^2\bar{a},e^{i\theta}z^2+a)
$$
where $c\in {\bf R}$, $\theta\neq \pi$ and $a$ a non-zero complex
number.

\bigskip

$(ii):$ $n=1$. We have the following

\bigskip

\textbf{Proposition 4.} \emph{Every cyclic, properly discontinuous
action on the unit disc $\Delta$ satisfies condition (1)}.

\bigskip

\emph{Proof.}  Let us collect some basic facts about the M\"obius
group (see eg. $\S\, 3.1$ in [16]). First of all, we recall
$$
{\rm
Aut}(\Delta)=\left\{\gamma(z)=\frac{\bar{p}z+\bar{q}}{qz+p}:|p|^2-|q|^2=1\right\}.
$$
If $\gamma\in {\rm Aut}(\Delta)$ is not the identity, then it is
 called elliptic, parabolic or hyperbolic according as
 $(p+\bar{p})^2<4,\,=4,\,>4$, and the eigenvalues of the coefficient matrix
$$
A=\left(
\begin{array}{ll}
\bar{p} & \bar{q}\\
q & p
\end{array}
\right)
$$
are two distinct conjugate complex numbers, two equal real numbers
or two distinct real numbers respectively.
  In the first case the fixed points are
 connected by a reflection in the unit circle, and in the other two
 cases the fixed points lie on the unit circle. Since $\Gamma=\{\gamma^k:k\in {\bf
 Z}\}$ acts properly discontinuously, we must rule out the elliptic case.

 $(ii)_1:$ Hyperbolic case.  Let $\rho_1,\rho_2$ be
the fixed points of $\gamma$ and $\lambda_1,\lambda_2$ the (real)
eigenvalues of the associated coefficient matrix $A$. If we
introduce the matrix
$$
B=\left(
\begin{array}{ll}
1 & -\rho_1\\
1 & -\rho_2
\end{array}
\right),
$$
then using an appropriate scalar factor $\alpha$, we will obtain the
relation
$$
B^{-1}\left(
\begin{array}{ll}
\mu & 0\\
0 & 1
\end{array}
\right) B=\alpha A
$$
where $\mu=\lambda_1/\lambda_2>0$.
 Simple computations show
 $$
 \alpha^k A^k=B^{-1}\left(
\begin{array}{ll}
\mu^k & 0\\
0 & 1
\end{array}
\right) B=\frac1{\rho_1-\rho_2}\left(
\begin{array}{ll}
-\rho_2\mu^k+\rho_1 & \rho_1\rho_2(\mu^k-1)\\
1-\mu^k & \rho_1 \mu^k-\rho_2
\end{array}
\right),
 $$
so that
$$
\gamma^k(0)=\frac{\rho_1\rho_2(\mu^k-1)}{\rho_1 \mu^k-\rho_2}.
$$
Since $\rho_1,\rho_2$ lie on the unit circle, we have
$$
\sum_{k\in {\bf Z}}(1-|\gamma^k(0)|^2)=\sum_{k\in {\bf
Z}}\left(1-\left|\frac{1-\mu^k}{\mu^k-\bar{\rho}_1\rho_2}\right|^2\right)=\sum_{k\in
{\bf Z}}\frac{2(1-{\rm
Re\,}(\bar{\rho}_1\rho_2))\mu^k}{|\mu^k-\bar{\rho}_1\rho_2|^2}<\infty
$$
because $\mu\neq 1$.

$(ii)_2:$ Parabolic case. Let $\lambda$ be the real double
eigenvalue of $A$. We may choose a matrix
$$ T=\left(
\begin{array}{ll}
a & b\\
c & d
\end{array}
\right)\ \ \ \ \ (ad-bc=1).
$$
such that
$$
\lambda^{-1} A=T^{-1}\left(
\begin{array}{cc}
1 & \nu \\
0 & 1
\end{array}
\right)T
$$
for certain $\nu\neq 0$. By a direct computation, we obtain
$$
\lambda^{-k}A^k=T^{-1}\left(
\begin{array}{cc}
1 & k\nu \\
0 & 1
\end{array}
\right)T=\left(
\begin{array}{cc}
1+cd\nu k & d^2\nu k \\
-c^2\nu k & 1-cd\nu k
\end{array}
\right),
$$
in particular,
$$
\lambda^{-1}\left(
\begin{array}{ll}
\bar{p} & \bar{q}\\
q & p
\end{array}
\right)=\left(
\begin{array}{cc}
1+cd\nu  & d^2\nu  \\
-c^2\nu  & 1-cd\nu
\end{array}
\right)
$$
 so that
 $$
 {\rm Re\,}(cd\nu)=0 \ \ \ {\rm and\ \ \ } |c|=|d|.
 $$
Clearly, $d\neq 0$. Hence
$$
\sum_{k\in {\bf Z}}(1-|\gamma^k(0)|^2)=\sum_{k\in {\bf
Z}}\left(1-\left|\frac{d^2\nu k }{1-cd\nu
k}\right|^2\right)=\sum_{k\in {\bf
Z}}\frac1{1+|d|^4|\nu|^2k^2}<\infty.
$$

\section{Proof of Theorem 2}

\ \ \ \ \ 5.1. First of all, we will derive a useful inequality. For
two distinct points $a,a'$ in $M$, we define
$$
C_{a,a'}=\frac{\sup\{-g_M(z,a):z\in
\partial D_{a,a'}\}}{\inf\{-g_M(z,a'):z\in
\partial D_{a,a'}\}}
$$
where $D_{a,a'}$ is a relatively compact open subsets containing
$a,a'$. Consider the plurisubharmonic function on $M$ defined by
$g_M(z,a)$ on $D_{a,a'}$  and $\max\{g_M(z,a),C_{a,a'}g_M(z,a')\}$
otherwise. It follows immediately from the definition of $g_M$ that
\begin{equation}
g_M(z,a)\ge C_{a,a'}g_M(z,a')\ \ \ \ \ {\rm for\ \ \ } z\in
M\backslash D_{a,a'}.
\end{equation}

 5.2. Let $a\in M$.   As $\Gamma$ is discrete,
we infer from (2), (5) that the inequality
$$
\sum_{\gamma\in \Gamma}g_M(\gamma(z),a)>-\infty
$$
holds  for all $z\in M\backslash \Gamma_a$. Put
$u(z)=\sum_{\gamma\in \Gamma}g_M(\gamma(z),a)$. This
$\Gamma-$invariant function descends to a negative plurisubharmonic
function on $M/\Gamma$ with a logarithmic pole $b:=\pi(a)$ where
$\pi$ is the natural projection, we claim that
$g_{M/\Gamma}(\cdot,b)$ exists. Since $a$ is arbitrary, we claim
that $M/\Gamma$ is pluricomplex hyperbolic.

\bigskip

5.3. Next, we are going to construct bounded plurisubharmonic
functions on $M/\Gamma$ which is \emph{strictly} plurisubharmonic at
one point. For $b\in M/\Gamma$, we may choose a coordinate ball
${\bf B}(b,1)$ at $b$ so that the inequalities
$$
\log|z-b|+C_1\le g_{M/\Gamma}(z,b)\le \log|z-b|+C_2
$$
hold in a open neighborhood of the closure of ${\bf B}(b,1)$. Here
$C_1\le C_2$ are two constants depending on $b$ but not on $z$. Take
$0<r<1$ so small that $\log r\le 2(C_1-C_2)$. Put
$$
v(z)=\left\{
\begin{array}{ll}
\log|z-b| & |z-b|<r \\
\max\{\log|z-b|,2(g_{M/\Gamma}(z,b)-C_1)\} & r\le |z-b|\le 1\\
2(g_{M/\Gamma}(z,b)-C_1) & |z-b|>1.
\end{array}
\right.
$$
Note that on the sphere $\{z:|z-b|=1\}$ we have
$$
2(g_{M/\Gamma}(z,b)-C_1)\ge 0= \log|z-b|
$$
while
$$
2(g_{M/\Gamma}(z,b)-C_1)\le 2(\log|z-b|+C_2-C_1)\le \log|z-b|
$$
holds on the sphere $\{z:|z-b|=r\}$. Thus $v$ is a well-defined
plurisubharmonic function on $M/\Gamma$ and $\rho:=\exp(2v)$ is a
bounded plurisubharmonic function on $M/\Gamma$ satisfying
$\rho(z)=|z-b|^2$ on $B(b,r)$.

\bigskip

5.4. \emph{Steinness.} Fix  a plurisubharmonic exhaustion function
$\psi$ on $M$. For $t>0$ we put $M_t=\{z\in M:\psi(z)<t\}$ and
$V_t=\pi(M_t)$. Again every $V_t$ is locally Stein relative compact
domain in $M/\Gamma$ and we have $M/\Gamma=\bigcup_{t>0}V_t$.
According to the argument in 3.4, it suffices to show that for every
$t>0$ there is a continuous strictly plurisubharmonic function in a
neighborhood of $\overline{V_t}$. By 5.3, there exists finite
coordinate balls ${\bf B}(w_i,r_i),\,i=1,\cdots,m$ and bounded
plurisubharmonic functions $\rho_i$ such that $\{{\bf
B}(w_i,r_i)\}_{1\le i\le m}$ cover $\overline{V_t}$ and  equalities
$\rho_i(z)=|z-w_i|^2$ hold on ${\bf B}(w_i,r_i)$. Put
$\rho=\rho_1+\cdots+\rho_m$. Then the positive current
$\partial\bar{\partial}\rho $ will dominate a continuous positive
$(1,1)-$form in a neighborhood of $\overline{V}_t$. By a standard
regularization of the function $\rho$, we get a smooth strictly
plurisubharmonic function in a neighborhood of $\overline{V_t}$.

\bigskip

5.5. In view of Corollary 2, it is natural to ask the following

\bigskip

\emph{Question.} For a Stein manifold, is pluricomplex hyperbolicity
equivalent to the existence of a negative $C^\infty$ strictly
plurisubharmonic function?

\section{Proof of Theorem 3}

\ \ \ \ \,6.1. \emph{Proof of 1)}. It is easy to construct a
negative plurisubharmonic function $\rho$ so that $-\rho\asymp
\delta_M$. Let $R$ be the diameter of $M$. Fix $w\in M$ and put
$$
C_w=\frac{\sup_{\{\rho(z)=\rho(w)/2\}}|\log|z-w|/R|}{\inf_{\{\rho(z)=\rho(w)/2\}}(-\rho(z))}.
$$
We can define a negative plurisubharmonic function on $M$ which
equals to $\log|z-w|/R$ for $\rho(z)\le \rho(a)/2$ and
$\max\{\log|z-w|/R,C_w\rho(z)\}$ otherwise. Hence
\begin{equation}
g_M(z,w)\ge C_w\rho(z)\ \ \ \ \ {\rm for\ \ \ } \rho(z)>\rho(w)/2.
\end{equation}
Since
 $M$ is convex, $g_M$ is symmetric (cf. [8]).
Thus for $z\in M\backslash\Gamma_a$ we have
$$
-\sum_{\gamma\in \Gamma} g_M(\gamma(z),a)=-\sum_{\gamma\in \Gamma}
g_M(z,\gamma^{-1}(a))  = -\sum_{\gamma\in \Gamma}
g_M(\gamma^{-1}(a),z)<\infty
$$
according as (6). Combining with (5), the conclusion follows.

\bigskip

6.2. \emph{Proof of 2)}. Usually one denotes by  ${\cal T}_{g,n}$
 the Teichm\"uller space of Riemann surfaces with genus $g$
and $n$ punctures. It is a bounded Stein domain in ${\bf
C}^{3g-3+n}$. Recently, Krushkal [7] observed that the pluricomplex
Green function relates to the Teichm\"uller  distance  by
$$
g_{{\cal T}_{g,n}}(z,w)=\log\tanh d_T(z,w),
$$
in particular, $g_{{\cal T}_{g,n}}$ is symmetric. Thus for $z\in
M-\Gamma_a$ we have
\begin{eqnarray*}
 &-&\sum_{\gamma\in \Gamma} g_{{\cal T}_{g,n}}(\gamma(z),a) = -\sum_{\gamma\in \Gamma}
g_{{\cal T}_{g,n}}(z,\gamma^{-1}(a))  = -\sum_{\gamma\in \Gamma}
g_{{\cal T}_{g,n}}(\gamma^{-1}(a),z)\\
&=& \sum_{\gamma\in \Gamma}\log\tanh d_T(\gamma^{-1}(a),z)  \le C_z
\sum_{\gamma\in \Gamma}e^{-2d_T(\gamma^{-1}(a),z)}\\
&\le &C_ze^{2d_T(z,a)} \sum_{\gamma\in
\Gamma}e^{-2d_T(\gamma^{-1}(a),a)}<\infty.
\end{eqnarray*}
By (5), the assertion follows.

\section{Proof of Theorem 4}

\ \ \ \ 7.1. \emph{Proof of Theorem 4}. It is well-known that the
series
$$
\sum_{\gamma\in \Gamma}|j_\gamma(z)|^2
$$
is always convergent at every $z\in \Omega$, moreover, it is
convergent uniformly on compact subsets of $\Omega$  (cf. Lemma 5 of
Chapter 10 in [13]). Here $j_\gamma$ denotes the Jacobian of
$\gamma$. On the other hand, from relation
\begin{equation}
K_\Omega(z)=K_\Omega(\gamma(z))|j_\gamma(z)|^2
\end{equation}
we conclude that
$$
\sum_{\gamma\in \Gamma}\frac1{K_\Omega(\gamma(z))}
$$
is also convergent uniformly on compact subsets of $\Omega$. Suppose
$p\in L_\Gamma$, namely, there exists a compact subset $K$,
$\{z_k\}\subset \Omega$ with $z_k\rightarrow p$ and
$\{\gamma_k\}\subset \Gamma$ so that $\gamma_k(z_k)\in K$. Observe
that
$$
\lim_{k\rightarrow \infty}K_\Omega(z_k)=\lim_{k\rightarrow
\infty}K_\Omega(\gamma_k^{-1}(\gamma_k(z_k)))=\infty,
$$
hence  $p\in \Lambda_\Omega$. For the second assertion, it suffices
to verify that $\Lambda_\Omega$ is closed in $\partial \Omega$.
Suppose $\Lambda_\Omega\ni p_k\rightarrow p$. Since for every $k$
there exists a point $z_k\in \Omega$ satisfying $|z_k-p_k|<1/k$ and
$K_\Omega(z_k)>k$, we have $p\in \Lambda_\Omega$.

\bigskip

7.2. It seems worthwhile to mention the following

\bigskip

\textbf{Proposition 5.} \emph{Let $\Omega$ be a bounded domain in
${\bf C}^n$ and $\Gamma$ a countable group of automorphisms of
$\Omega$. Then $\Gamma$ is properly discontinuous iff there exists
one point $a\in \Omega$ such that}
$$
\sum_{\gamma\in \Gamma} |j_\gamma(a)|^2<\infty.
$$

\emph{Proof.} It suffices to verify the "if" part. By (7) we have
$$
\sum_{\gamma\in \Gamma}\frac1{K_\Omega(\gamma(a))}<\infty.
$$
As $1/K_\Omega$ is a positive continuous function, the assertion
follows immediately from Proposition 1.

\section{Proofs of Theorem 5, 6, 7}

\ \ \ \ 8.1. \emph{Proof of Theorem 5.} By Richberg's theorem [11],
we may assume $\psi$ is negative, $C^\infty$ strictly
plurisubharmonic on $\Omega$. Put
$$
u(z)=\sup_{\gamma\in \Gamma} \psi(\gamma(z)),\ \ \ z\in \Omega.
$$
Clearly, $u$ is $\Gamma-$invariant. Since $\psi(z)\rightarrow
-\infty$ as $z\rightarrow \cup_{a\in D}{\rm C}\Gamma_a$, we conclude
that for any open subset $W\subset\subset \Omega$ there exists a
finite subset $\Gamma'$ of $\Gamma$ so that
$$
\psi(\gamma(z))<\max_{\gamma\in \Gamma'} \psi(\gamma(z)),\ \ \
\forall\, z\in W,\,\gamma\in \Gamma\backslash\Gamma'.
$$
Thus $u(z)=\max_{\gamma\in \Gamma'} \psi(\gamma(z))$ for $z\in W$ so
that it is continuous and strictly plurisubharmonic there. As $W$ is
arbitrarily chosen, $u$ descends to a negative, continuous and
strictly plurisubharmonic function on $\Omega/\Gamma$. The rest
argument is same as 3.4.

\bigskip

8.2. \emph{Proof of Theorem 6.}  It is well-known that  elements of
${\rm Aut}({\bf B})\backslash\{I\}$ are divided into three classes:
those have fixed points in ${\bf B}$ are called \emph{elliptic},
those have one or two fixed points on $\partial {\bf B}$ and no
fixed points in ${\bf B}$ are called \emph{parabolic} or
\emph{hyperbolic} respectively (see e.g., [12]). Let $\Gamma$ be an
infinite, cyclic, free and properly discontinuous subgroup of ${\rm
Aut}({\bf B})$ and $\gamma$ the generator of $\Gamma$. Clearly,
$\gamma$ can not be elliptic.

By Proposition 2, it suffices to show that ${\rm C}\Gamma_0$
contains at most two elements.  Let $p\in {\rm C}\Gamma_0$, namely,
there there is a subsequence $\gamma^{k_j}(0)\rightarrow p$ as
$j\rightarrow \infty$. Thus we have
$$
\gamma^{k_j+1}(0)=\gamma(\gamma^{k_j}(0))\rightarrow \gamma(p)
$$ since $\gamma$ is continuous on the closure of $\bf B$.
On the other hand, Proposition 3 implies
$$
\gamma^{k_j+1}(0)=\gamma^{k_j}(\gamma(0))\rightarrow p.
$$
Hence $\gamma(p)=p$, i.e., $p$ is a fixed point of $\gamma$ on
$\partial {\bf B}$.  The proof is complete.

\bigskip

8.3. \emph{Proof of Theorem 7.} Every element in ${\rm
Aut\,}(\Delta^2)$ is of form
$$
\gamma(z^1,z^2)=(\gamma_1(z^1),\gamma_2(z^2))\ \ \ {\rm or}\ \ \
\gamma(z^1,z^2)=(\gamma_2(z^2),\gamma_1(z^1))
$$
where $\gamma_1,\gamma_2\in {\rm Aut\,}(\Delta)$. In the first case,
we have
$$
\gamma^k(z^1,z^2)=(\gamma_1^k(z^1),\gamma_2^k(z^2)).
$$
As $\Gamma=\{\gamma^k:k\in {\bf Z}\}$ acts properly discontinuously,
we conclude that there exists at least one component of $\gamma$,
say $\gamma_1$, which can not be elliptic. The argument in 8.2
yields
$$
{\rm C}\Gamma_z\subset {\rm Fix}(\gamma_1)\times \Delta,\ \ \
\forall\, z\in \Delta^2
$$
where ${\rm Fix}(\gamma_1)$ is the set of fixed points $\gamma_1$.
Hence the condition of Theorem 5 is verified (e.g., take
$\psi(z)=\sum_{p\in {\rm Fix}(\gamma_1)}\log|z^1-p|$).

In the second case, we have
$$
\gamma^{2m}(z^1,z^2)=((\gamma_2\circ\gamma_1)^m(z^1),(\gamma_1\circ\gamma_2)^m(z^2)).
$$
For every $z\in {\Delta}^2$, the set ${\rm C}\Gamma_z$ consists of
the cluster points of $\{\gamma^{2m}(z):m\in {\bf Z}\}$ and
$\{\gamma^{2m+1}(z):m\in {\bf Z}\}$. From relation
$\gamma^{2m+1}=\gamma\circ \gamma^{2m}$, we see that that $\gamma$
maps the set of cluster points of the first to the second. On the
other hand, since
\begin{eqnarray*}
p\in {\rm Fix}(\gamma_2\circ \gamma_1) & \Rightarrow &
\gamma_1(p)\in {\rm Fix}(\gamma_1\circ\gamma_2)\\
p\in {\rm Fix}(\gamma_1\circ \gamma_2) & \Rightarrow &
\gamma_2(p)\in {\rm Fix}(\gamma_2\circ\gamma_1),
\end{eqnarray*}
neither $\gamma_1\circ\gamma_2$ nor $\gamma_2\circ\gamma_1$ is
elliptic. The argument in 8.2 shows that the cluster points of
$\{\gamma^{2m}(z):m\in {\bf Z}\}$ are contained in ${\rm
Fix}(\gamma_2\circ\gamma_1)\times {\rm Fix}(\gamma_1\circ\gamma_2)$
which has at most four elements, hence $\cup_{z\in {\Delta}^2}{\rm
C}\Gamma_z$ has at most eight elements. The proof is complete.

\bigskip

8.4. Finally, let us present another example. Let
$$
\Omega=\{(z^1,z')\in {\bf C}\times {\bf
C}^{n-1}:|z^1|^2+\psi(z')<1\}
$$
where $\psi$ enjoys the following properties: i) $\psi$ is
$C^\infty$ and plurisubharmonic in ${\bf C}^{n-1}$; ii) there exist
positive numbers $\alpha,\,C$ so that $\psi(z')\ge C|z'|^\alpha$;
iii) $\psi$ is weighted homogeneous, i.e., there are positive
integers $\tau_1,\tau_2,\cdots,\tau_{n}$ such that
$$
\psi(c^{\tau_2}z^2,\cdots,c^{\tau_n}z^n)=|c|^{\tau_1}\psi(z'),\ \ \
\forall\,c\in {\bf C}.
$$
Property i), ii) implies that $\Omega$ is a bounded Stein domain.
Consider the Cayley transformation
$$
\Phi:\Omega\rightarrow {\bf H}:=\{z\in {\bf C}^n:{\rm
Im\,}z^1>\psi(z')\}
$$
$$
z\rightarrow
\left(i\frac{1+z^1}{1-z^1},\frac{iz^2}{(1-z^1)^{2\tau_2/\tau_1}},\cdots,\frac{iz^n}{(1-z^1)^{2\tau_n/\tau_1}}\right).
$$
Simple computations show
$$
\Phi^{-1}(z)=\left(\frac{z^1-i}{z^1+i},\frac{2^{2\tau_2/\tau_1}i^{2\tau_2/\tau_1-1}}{(z^1+i)^{2\tau_2/\tau_1}}z^2,\cdots,
\frac{2^{2\tau_n/\tau_1}i^{2\tau_n/\tau_1-1}}{(z^1+i)^{2\tau_n/\tau_1}}z^n\right).
$$
Fix a positive number $r\neq 1$. Put
$$
\tilde{\gamma}(z)=(r^{\tau_1}
z^1,r^{\tau_2}z^2,\cdots,r^{\tau_n}z^n)\in {\rm Aut\,}({\bf H})
$$
and $\Gamma=\{\gamma^k: k\in {\bf Z}\}$ where $\gamma=\Phi^{-1}\circ
\tilde{\gamma}\circ \Phi$. It is not difficult to see ${\rm
C}\Gamma_0=\{(1,0'),(-1,0')\}$ and the function $|z^1|^2-1$ provides
a barrier at $(1,0')$ and $(-1,0')$. Thus Theorem 5 applies.

\end{document}